\documentstyle[12pt,amsfonts]{article}
\def\N{{\Bbb N}}

\def\Z{{\Bbb Z}}

\begin{document}
\baselineskip=22pt
\title{\Large\bf The Zeta function, Periodic Points and Entropies\\
of the Motzkin Shift}
\author{Kokoro Inoue\\
{\small \em Yakuin 2-3-26, Chuo-ku, Fukuoka 810-0022, Japan}\\
{\small \em kokoro$\_{}$inoue@gakushikai.jp}}
\date{}
\maketitle
\vspace{1cm}

\begin{abstract}
We discuss a method of calculating the zeta function of subshifts which have a presentation by a finite directed graph labeled by elements of the associated inverse semigroup. This class of subshifts is introduced as a class of property A subshifts(T.Hamachi, K.Inoue and W.Krieger/Subsystems of finite type and semigroup invariants of subshifts, preprint), and the Dyck shift and the Motzkin shift are representative subshifts in this class. The exact number of the periodic points and entropies are also given by this method, and these values are used in the embedding condition for an irreducible shift of finite type into a subshift in this class.
\end{abstract}

\newtheorem{definition}{Definition}[section]
\newtheorem{lemma}{Lemma}[section]
\newtheorem{theorem}{Theorem}[section]
\newtheorem{pr}{Proposition}[section]
\newtheorem{Motivation}{Motivation}[section]
\newtheorem{cor}{Corollary}[section]
\newtheorem{rem}{Remark}[section]


\section{Introduction}
~
Put $\Sigma=\{\lambda_1,\lambda_2,\cdots,\lambda_M,
\rho_1,\rho_2,\cdots,\rho_M,1_1,1_2,\cdots,1_N\},\ M,N\in\N$. Let ${\cal M}(M,N)$ be a monoid (with zero) with generators 
$\lambda_i,\rho_i,\ 1\leq i\leq M,\ 1_j,\ 1\leq j\leq N$ and ${\bf 1}$.
The relations are:

$$\begin{array}{ll}
\lambda_i\cdot\rho_i={\bf 1},&1\leq i\leq M,\\
1_i\cdot1_j={\bf 1},&1\leq i,j\leq N,\\
\lambda_i\cdot\rho_j=0,&1\leq i,j\leq M,\ i\ne j,\\
{\bf 1}\cdot\alpha=\alpha\cdot{\bf 1}=\alpha,&\alpha\in\Sigma\cup\{{\bf 1}\},\\
1_i\cdot\alpha=\alpha\cdot1_i=\alpha,&1\leq i\leq N,\ \alpha\in\Sigma\cup\{{\bf 1}\},\\
0\cdot\alpha=\alpha\cdot0=0,&\alpha\in\Sigma\cup\{{\bf 1}\},\\
0\cdot0=0.
\end{array}$$

Set a mapping ${\it red}:\Sigma^\ast\longrightarrow{\cal M}(M)$ by
$${\it red}(\alpha)=\alpha_1\cdot\alpha_2\cdot\dots\cdot\alpha_n\ \ \mbox{for}\ \ \alpha=\alpha_1\alpha_2\dots\alpha_n,\ n\in\N$$
and 
$${\it red}(\varepsilon)={\bf 1}\ \ \mbox{where}\ \varepsilon\ \mbox{is\ the\ empty\ word.}$$
The Motzkin shift ${\bf M}(M,N)$ is defined by
$${\bf M}(M,N)=\{x\in\Sigma^{\Z}:{\it red}(x_ix_{i+1}\cdots x_j)\ne0\ \mbox{for}\ \forall i\leq\forall j\}.$$

Therefore we can regard the Motzkin shift as a shift defined by a simple directed graph $G$ which has one vertex and $(2M+N)$-loops named by the elements of the set $\Sigma$, and the loop named $\lambda_i,\ \rho_i(1\leq i\leq M)$ carry the labels $\lambda,\ \rho$ respectively, the loop named $1_i(1\leq i\leq N)$ carries the label {\bf 1}, that is, ${\cal M}(M)$ is the Dyck type inverse monoid. And such a presentation of subshifts is called ${\it S}$-presentation in [HIK], where ${\it S}$ is an inverse semigroup of Dyck type. Note that if $N=0$, the monoid is the Dyck monoid ${\cal D}_M$ and the subshift ${\bf M}(M,0)$ is the Dyck shift ${\bf D}_M$ in [HI].

In [HIK], a necessary and sufficient condition for embedding of an irreducible subshift of finite type into a property A subshift that allows an ${\it S}$-presentation is obtained. This condition is extension of the embedding condition whose target is the Dyck shift in [HI], and which consists of periodic condition and the entropy condition. For the Dyck shift ${\bf D}_M$ the entropy, $h({\bf D}_M)=\log(M+1)$, is know by W.Krieger [Kr] and the number of $n$-periodic points ($n\in\N$) is obtained in [I] and [HI]. As known, the zeta function of the Dyck shift is calculated by circular codes [Ke]. Also from the number of periodic points, the entropy and the zeta function are computed directly[I].

In section 2, we consider the zeta function(Proposition 2.2), the number of periodic points(Proposition 2.3) and the entropy(Proposition 2.4) of the Motzkin shift. Here it is a crucial point that periodic points are classified by images of periodic defining blocks by the mapping ${\it red}$. 

Put $X$ a subshift over $\Sigma$ which allows an ${\it S}$-presentation. If $\alpha=\alpha_1\alpha_2\cdots\alpha_n\in {\cal B}_n(X)$, $n\geq1$ satisfies $\alpha^\infty\in X$, we call $\alpha^\infty$ an $n$-periodic point of $X$ and $\alpha$ is called an $n$-periodic defining block of $\alpha^\infty$ where ${\cal B}_n(X)$ denotes a set of length $n$ blocks and ${\cal B}(X)=\bigcup_{n=1}^{\infty}{\cal B}_n(X)$. Any block of the Motzkin shift $\alpha$ is given as a product $\alpha=\alpha_+\alpha_-$ such that ${\it red}(\alpha_+)\in{\cal M}^+(M)$ and ${\it red}(\alpha_-)\in{\cal M}^-(M)$, where ${\cal M}^+(M)$ is the free monoid generated by $\rho_i,\ 1\leq i \leq M$ and ${\cal M}^-(M)$ is the free monoid generated by $\lambda_i,\ 1\leq i \leq M$. For a periodic defining block $\alpha=\alpha_+\alpha_-$ we call ${\it red}(\alpha_-\alpha_+)$ a multiplier of $\alpha$ or of $\alpha^\infty$. For a block $\alpha=\alpha_+\alpha_-$, if a multiplier of $\alpha$ is either positive $({\it i.e.}\ {\it red}(\alpha_-\alpha_+)\in{\cal M}^+)$, negative $({\it i.e.}\ {\it red}(\alpha_-\alpha_+)\in{\cal M}^-)$ or neutral $({\it i.e.}\ {\it red}(\alpha_-\alpha_+)={\bf 1})$, then $\alpha$ is a periodic defining block, and vice versa. In this paper, $P_n(X)$ denotes the set of $n$-periodic points and $p_n(X)=\sharp P_n(X)$, $n\geq1$.

Using these property of the Motzkin shift (the Dyck shift), bijections are constructed between a set of periodic points belonging to some class of this classification and an edge shift. The bijection gives the generating function of a circular code that is needed in the calculation of the zeta function. This generating function can also be obtained by using the Sch$\ddot{\rm u}$tzenberger's method.

The entropies $h^+$ and $h^-$ are introduced for the embedding condition in [HIK]. In fact, the Dyck shift and the Motzkin shift satisfy $h^+=h^-$ because of their symmetry. In section 3, we show a simple example satisfying $h^+>h^-$.

At the time when the author reached the results, she had not knowledge of the \\Sch${\rm \ddot{u}}$tzenberger's method. She would like to thank W.Krieger, who let me know the \\
Sch$\ddot{\rm u}$tzenberger's method.

\bigskip

\section{The zeta function of the Motzkin shift}
~
The Motzkin shift ${\bf M}(M,N)$ is defined by a directed graph $G=G(M,N)$ with one vertex and $(2M+N)$-loops named by $\lambda_i,\rho_i,\ 1\leq i\leq M,\ 1_j,\ 1\leq j\leq N$. Put $G'=G'(M,N)$ a graph which is gotten by removing the loops named $\lambda_2,\lambda_3,\cdots,\lambda_M$ from the graph $G$ and $X_{G'}$ the edge shift of $G'$. Then we can construct a bijection as follows.

\begin{pr}
For $n\geq1$ there exists bijections
$$
\left\{x\in P_n({\bf M}(M,N)):
\begin{array}{l}
\mbox{a multiplier of}\ x\ \mbox{is either positive, neutral or}\\
\mbox{a positive power of}\ \lambda_1
\end{array}
\right\}
\longrightarrow {\cal B}_n(X_{G'}).
$$
\end{pr}

\noindent{\bf Proof}\quad Take $x=\alpha^\infty\in P_n({\bf M}(M,N))$, $\alpha=\alpha_+\alpha_-$ such that a multiplier of $\alpha$ is positive or neutral or a positive power $\lambda_1$. Replacing all $\lambda_i,\ 2\leq i\leq M$ in $\alpha$ with $\lambda_1$, we get a block of $X_{\cal G'}$, and obviously it gives a bijection. (See Proposition 2.4 in [HI],[I].)
\hfill$\Box$

\bigskip

\noindent Set\par
$R=\{\rho_1,\rho_2,\cdots,\rho_M,1_1,1_2,\cdots,1_N\}\cup\{\varepsilon\},$\par
$L=\{\lambda_1,\lambda_2,\cdots.\lambda_M,1_1,1_2,\cdots,1_N\}\cup\{\varepsilon\}\quad \mbox{and}\quad \widetilde L=\{\lambda_1,1_1,1_2,\cdots,1_N\}\cup\{\varepsilon\}.$\par
$n\geq1,$\par
$E_n=\{w\in {\cal B}_n({\bf M}(M,N)):w=\lambda_iv\rho_i,\ 1\leq i\leq M, \mbox{where}\ {\it red}(v)=\mbox{\bf 1}\}$\par
and $E=\bigcup_{n=1}^{\infty}E_n$.\par
$E_+=\{wu\in {\cal B}({\bf M}(M,N)):w\in E, u\in R^\ast\},$\par
$E_-=\{uw\in {\cal B}({\bf M}(M,N)):w\in E, u\in L^\ast\}$ and\par
$\widetilde E_-=\{uw\in {\cal B}({\bf M}(M,N)):w\in E, u\in \widetilde L^\ast\}.$\par
$F=\{wv\in {\cal B}({\bf M}(M,N)):w\in E, v\in \{1_1,1_2,\cdots,1_N\}^\ast\}.$\\
\bigskip

Then these are circular codes. And particularly codes $E$ are called the Motzkin codes(the Dyck codes for $N=0$). In general, for a circular code $C$, the zeta function of $C^{\infty}$, $\zeta(C^{\infty},z)$ is given by the generating function of $C$, which is defined by $f(C,z)=\sum_{u\in C}z^{|u|}$, as follows [St]
$$\zeta(C^{\infty},z)=\frac1{1-f(C,z)}.$$

In [Ke], the zeta function of the Dyck shift is calculated by using circular codes $E$, $E_+$ and $R$. In the same way, the zeta function of the Motzkin shift is gotten by the circular codes defined above. But, in fact, the true nature of the calculation of these zeta functions is to calculate the generating function of the Motzkin code(the Dyck code) $E$.

\begin{theorem}
The generating function of the Motzkin code $f(E,z)$ satisfies
$$f(E,z)^2+(Nz-1)f(E,z)+Mz^2=0.$$
Hence $f(E,z)=\frac12\{1-Nz-\sqrt{(1-Nz)^2-4Mz^2}\}.$
\end{theorem}

\noindent{\bf Proof}\quad From Proposition 2.1 and $P_n(E_+^\infty)\cap P_n(\widetilde E_-^\infty)=P_n(F^\infty)$ and $P_n(R^\infty)\cap P_n(\widetilde L^\infty)=\{1_1,1_2,\cdots,1_N\}^n$, for $n\geq1$
\begin{equation}
p_n(X_{G'})=p_n(E_+^\infty)+p_n(\widetilde E_-^\infty)-p_n(F^\infty)+p_n(R^\infty)+p_n(\widetilde L^\infty)-N^n.
\end{equation}
Since $E_+, \widetilde E_-, F, R$ and $\widetilde L$ are circular codes, each zeta function is given by each generating function and the generating functions are $f(R,z)=(M+N)z$, $f(\widetilde L,z)=(N+1)z$ and
\begin{center}
$f(E_+,z)=f(E,z)\sum_{k=0}^\infty |R|^kz^k=f(E,z)/\{1-(M+N)z\}$.
\end{center}
Similarly, $f(\widetilde E_-,z)=f(E,z)/\{1-(N+1)z\}$ and $f(F,z)=f(E,z)/(1-Nz)$.
Therefore from (1)
\begin{eqnarray*}
\frac1{1-(1+M+N)z}
&=&\frac{\zeta(E_+^\infty,z)\zeta(\widetilde E_-^\infty,z)\zeta(R^\infty,z)\zeta(\widetilde L^\infty,z)(1-Nz)}{\zeta(F^\infty,z)}\\
&=&\frac{1-Nz-f(E,z)}{\{1-(M+N)z-f(E,z)\}\{1-(N+1)z-f(E,z)\}}.
\end{eqnarray*}
Thus the proof is complete.\hfill$\Box$

\begin{pr}
The zeta function of the Motzkin shift is
$$\zeta({\bf M}(M,N),z)=\frac{2\{1-Nz+\sqrt{(1-Nz)^2-4Mz^2}\}}{\{1-(2M+N)z+\sqrt{(1-Nz)^2-4Mz^2}\}^2}.$$
\end{pr}

\noindent{\bf Proof}\quad The number of periodic points of the Motzkin shift is \begin{equation}
p_n({\bf M}(M,N))=p_n(E_+^\infty)+p_n(E_-^\infty)-p_n(F^\infty)+p_n(R^\infty)+p_n(L^\infty)-N^n.
\end{equation}
By the symmetry of ${\bf M}(M,N)$,
$$\zeta({\bf M}(M,N),z)
=\frac{\zeta(E_+^\infty,z)^2\zeta(R^\infty,z)^2(1-Nz)}{\zeta(F^\infty,z)}\\
=\frac{(1-f(F,z))(1-Nz)}{(1-f(E_+,z))^2(1-f(R,z))^2}.$$
From theorem 2.1, the zeta function is obtained.\hfill$\Box$

\bigskip

In [I] and [HI], the number of periodic points of the Dyck shift ${\bf D}_M$ is
\begin{equation}
p_n(D_M)=\left\{ 
\begin{array}[c]{lll}
2\{(M+1)^n-\sum_{i=0}^{n/2}{n \choose i}M^i\}+{n \choose \frac{n}2}M^{\frac{n}2}&\quad&\mbox{if}\ n\ \mbox{is\ even},\\
&&\\
2\{(M+1)^n-\sum_{i=0}^{(n-1)/2}{n \choose i}M^i\}&\quad&\mbox{if}\ n\ \mbox{is\ odd}.
\end{array}\right.
\end{equation}
For the Motzkin shift, from (1),(2) and symmetry {\it i.e.} $p_n(E_-^\infty)=p_n(E_+^\infty)$ and $p_n(L^\infty)=p_n(R^\infty)$, $n\geq1$, the following is established:\\
\begin{equation}
p_n({\bf M}(M,N))=2\{p_n(X_{G'})-p_n(\widetilde E_-^\infty)-p_n(\widetilde L^\infty)\}+p_n(F^\infty)+N^n,\ n\geq1.
\end{equation}
Notice that $p_n(F^\infty)=p_n(E^\infty)$, $n\geq1$ for the Dyck shift ({\it i.e.}$N=0$). To compute the number of $n$-periodic points($n\geq1$). it is sufficient to compute $p_n(\widetilde E_-^\infty)+p_n(\widetilde L^\infty)$ and $p_n(E^\infty)$.(See the proof of Theorem 2.5 in [HI],[I].) For the Motzkin shift, it is possible to compute these values directly, but we can also obtain $p_n({\bf M}(M,N)),\ n\geq1$, easily by using (3) as follows.

\begin{pr}\quad For $n\geq1$
\begin{eqnarray*}
\lefteqn{p_n({\bf M}(M,N))}\\
&&=\left\{
\begin{array}[c]{ll}
2\{(M+N+1)^n-\sum_{i=0}^{n/2}{n \choose 2i}N^{2i}\sum_{j=0}^{n/2-i}{n-2i \choose j}M^j&\\
-\sum_{i=0}^{n/2-1}{n \choose 2i+1}N^{2i+1}\sum_{j=0}^{n/2-i-1}{n-2i-1 \choose j}M^j\}&\\
+\sum_{i=0}^{n/2}{n \choose n-2i}{2i \choose i}M^iN^{n-2i}&\mbox{if}\ n\ \mbox{is\ even},\\
&\\
2\{(M+N+1)^n-\sum_{i=0}^{n/2}{n \choose 2i}N^{2i}\sum_{j=0}^{(n-1)/2-i}{n-2i \choose j}M^j&\\
-\sum_{i=0}^{(n-1)/2}{n \choose 2i+1}N^{2i+1}\sum_{j=0}^{(n-1)/2-i}{n-2i-1 \choose j}M^j\}&\\
+\sum_{i=0}^{(n-1)/2}{n \choose n-2i}N^{n-2i}{2i \choose i}M^i&\mbox{if}\ n\ \mbox{is\ odd}.\end{array}\right.
\end{eqnarray*}
\end{pr}

\noindent{\bf Proof}\quad The number of periodic points of ${\bf M}(M,N)$ is given by the number of periodic points of the Dyck shift $p_n(D_M)$ as follows. 
$$p_n({\bf M}(M,N))=\sum_{i=0}^{n-1}{n \choose i}N^ip_{n-i}(D_M)+N^n.$$
Thus the proof is complete.\hfill$\Box$

\begin{pr}\quad The entropy of the Motzkin shift ${\bf M}(M,N)$ is 
$$h({\bf M}(M,N))=\log(M+N+1).$$
\end{pr}

\noindent{\bf Proof}\quad Since $\widetilde E_-\supset F$ and $p_n(\widetilde L)=(N+1)^n,\ n\geq1$, from (4)
$$2p_n(X_{G'})\geq p_n({\bf M}(M,N)),\ n\geq1.$$
Since $E_-\supset\widetilde E_-$ and $L\supset\widetilde L$,\ $n\geq1$, from (1),(2)
$$p_n({\bf M}(M,N))\geq p_n(X_{G'}),\ n\geq1.$$
Then $h({\bf M}(M,N))=h(X_{G'})=\log(M+N+1).$\hfill$\Box$

\bigskip

\section{A subshift whose entropies $h^+$ and $h^-$ differ}
~
The necessary and sufficient condition for embedding an irreducible subshifts of finite type into a property A subshift $X$ that allows an ${\it S}$-presentation in [HIK](Theorem 5.10) consists of periodic conditions and entropy conditions. And entropies $h^+(X)$, $h^-(X)$ are introduced for the entropy condition. They are defined by
$$h^+(X)=\liminf_{n\to\infty}\frac1n\log p_n^+(X),$$
where $p_n^+(X)$ is the number of $n$-periodic points of $X$ whose periodic defining block have a positive or neutral multiplier[HIK]. The entropy $h^-(X)$ is defined alike. The Dyck shift and the Motzkin shift satisfy $h^+=h^-$ because of their symmetry. Here we show an example possessing $h^+>h^-$.

Put $\Sigma_1=\{\lambda,\xi,\rho_1,\rho_2,\cdots,\rho_M,\eta_1,\eta_2,\cdots,\eta_M\},\ M\in\N$ and ${\cal M}_1(M)$ a monoid (with zero) with generators $\lambda,\xi,\rho_i,\eta_i, 1\leq i\leq M$ and {\bf 1} defined by
$$\lambda\cdot\rho_i={\bf 1},\quad \lambda\cdot\eta_i=0,\quad \xi\cdot\eta_i={\bf 1},\quad \xi\cdot\rho_i=0\ (1\leq i\leq M),$$
$${\bf 1}\cdot\alpha=\alpha\cdot{\bf 1}=\alpha,\quad 0\cdot\alpha=\alpha\cdot0=0\quad \mbox{and}\quad 0\cdot0=0.$$

Therefore monoids ${\cal M}_1^+(M)$ and ${\cal M}_1^-(M)$ are generated by $\rho_i,\eta_i$, $1\leq i\leq M$, {\bf 1} and $\lambda,\xi$, {\bf 1} respectively according to the above relations.

Set $X(M)=\{x\in\Sigma_1^{\Z}:{\it red}(x_ix_{i+1}\cdots x_j)\ne0\ \mbox{for}\ \forall i\leq\forall j\}.$ If $M=1$, $X(1)$ is the Dyck shift $D_2$ itself. This subshift $X(M)$ has an underlying graph $G_1=G_1(M)$ which has one vertex and 2(M+1)-loops named by the elements of $\Sigma_1$, and the loops with the names $\lambda$, $\xi$, $\rho_i$ and $\eta_i\ (1\leq i\leq M)$ carry the labels $\lambda$, $\xi$, $\rho$ and $\eta$ respectively.

\bigskip

\begin{theorem}
$$h^+(X(M))=\log(2M+1)\quad \mbox{and}\quad h^-(X(M))=\log(M+2).$$
Hence $h^+(X)>h^-(X)$ for $M\geq2$.
\end{theorem}

For proving Theorem 3.1, we prepare two graphs $G'_1$ and $G''_1$. A graph $G'_1$ has one vertex and (2M+1)-loops with the names $\lambda$, $\rho_i$, $\eta_i$ ($1\leq i\leq M$), and put the edge shift $X_{G'_1}$. A graph $G''_1$ has one vertex and (M+2)-loops with the names $\lambda$, $\xi$, $\rho_i$ ($1\leq i\leq M$), and put the edge shift $X_{G''_1}$. We can construct two bijections between some sets of $n$-periodic points and ${\cal B}_n(X_{G'_1})$ or ${\cal B}_n(X_{G''_1})$ respectively.

\begin{lemma}
For $n\geq1$ there exists bijections
\begin{equation}
\left\{x\in P_n(X(M)):
\begin{array}{l}
\mbox{a multiplier of}\ x\ \mbox{is positive, neutral or}\\
\mbox{a positive power of}\ \lambda
\end{array}
\right\}
\longrightarrow {\cal B}_n(X_{G'_1})
\end{equation}
and
\begin{equation}
\left\{x\in P_n(X(M)):
\begin{array}{l}
\mbox{a multiplier of}\ x\ \mbox{is negative, neutral or}\\
\mbox{a positive power of}\ \rho_1,\rho_2,\cdots,\rho_M
\end{array}
\right\}
\longrightarrow {\cal B}_n(X_{G''_1}).
\end{equation}
\end{lemma}

\noindent{\bf Proof}\quad The proof is done in the same way as of Proposition 2.1. \hfill$\Box$

\bigskip

\noindent{\bf Proof of Theorem 3.1}\quad Put\\
$R_1=\{\rho_1,\rho_2,\cdots,\rho_M,\eta_1,\eta_2,\cdots,\eta_M\}\cup\{\varepsilon\}$,\quad $L_1=\{\lambda,\xi\}\cup\{\varepsilon\}\ \mbox{and}\ \widetilde L_1=\{\lambda\}\cup\{\varepsilon\}.$\\
And circular codes $E, E_\pm, \widetilde E_-$ and $F$ are defined by the $R_1$ and $L_1$ in the same way as those of the Motzkin shift in the previous section.\par
From Lemma 3.1 (5), for $n\geq1$
\begin{eqnarray}
p_n(X_{G'_1})&=&p_n(E_+^\infty)+p_n(\widetilde E_-^\infty)-p_n(E^\infty)+p_n(R_1^\infty)+p_n(\widetilde L_1^\infty)\\
&>&p_n^+(X(M))\nonumber,
\end{eqnarray}
since $p_n^+(X(M))=p_n(E_+^\infty)+p_n(R_1^\infty)$, $n\geq1$ and $\widetilde E_-\supset E$.\par
On the other hand, since $p_n(E_+^\infty)>p_n(\widetilde E_-^\infty)$,\quad $p_n(R_1^\infty)>p_n(\widetilde L_1^\infty)$, $n\geq1$ and (7),
$$2p_n^+(X(M))>p_n(X_{G'_1}).$$
Then $h^+(X(M))=h(X_{G'_1})=\log(2M+1)$.\par
Set $\widetilde R_1=\{\rho_1,\rho_2,\cdots,\rho_M\}\cup\{\varepsilon\}$ and $\widetilde E_+=\{wu\in B(X(M)):w\in E, u\in \widetilde R_1^\ast\}$. From the bijection (6)
\begin{equation}
p_n(X_{G''_1})=p_n(\widetilde E_+^\infty)+p_n(E_-^\infty)-p_n(E^\infty)+p_n(\widetilde R_1^\infty)+p_n(L_1^\infty).
\end{equation}
Similarly, for $p_n^-(X(M))=p_n(E_-^\infty)+p_n(L_1^\infty)$, we obtain $p_n(X_{G''_1})\geq p_n^-(X(M))\geq\frac12 p_n(X_{G''_1})$. Then $h(X(M))=h(X_{G''_1})=\log(M+2)$.\hfill$\Box$

\bigskip

The zeta function $\zeta(X(M))$ and the number of periodic points $p_n(X(M))$ are also computable.

\begin{lemma}
The generating function of $E$, $f(E,z)$ is identical with of the Dyck shift ${\bf D}_{2M}$, that is, $f(E,z)=\displaystyle{\frac12}(1-\sqrt{1-8Mz^2})$.
\end{lemma}

\noindent{\bf Proof}\quad The proof is done in the same way as Proposition 2.1 by using the bijection (5) in Lemma 3.1. Here, notice that the graph $G'_1$ defining ${\cal B}_n(X_{G'_1})$ is identical with the graph $G'(2M,0)$. Therefore Lemma 3.1 is immediately established from Proposition 2.2. \hfill$\Box$

\begin{pr}
The zeta function of the shift $X(M)$ is
$$\zeta(X(M),z)=\frac{2(1+\sqrt{1-8Mz^2)}}{(1-4Mz+\sqrt{1-8Mz^2})(1-4z+\sqrt{1-8Mz^2})}.$$
\end{pr}

\noindent{\bf Proof}\quad The proof is done in the same way as of Proposition 2.3, using the circular codes defined in the proof of Theorem 3.1 and the generating function $f(E,z)$ obtained in Lemma 3.2.\hfill$\Box$

\bigskip

The zeta functions for $p_n^+(X(M))$ and $p_n^-(X(M))$ are also calculated. Set $$\zeta^\pm(X(M),z)=\exp\sum_{n=1}^{\infty}\frac{p_n^\pm(X(M))}nz^n.$$

\begin{pr}
$$\zeta^+(X(M),z)=\frac2{1-4Mz+\sqrt{1-8Mz^2}}$$
and
$$\zeta^-(X(M),z)=\frac2{1-4z+\sqrt{1-8Mz^2}}.$$
\end{pr}

\noindent{\bf Proof}\quad Since $p_n^+(X(M))=p_n(E_+^\infty)+p_n(R_1^\infty)$, then $\zeta^+(X(M),z)=\zeta(E_+^\infty,z)\zeta(R_1^\infty,z)$. From the generating functions of $E_+$ and $R_1$, the zeta function $\zeta^+(X(M),z)$ is obtained.\par
In the same way, from $p_n^-(X(M))=p_n(E_-^\infty)+p_n(L_1^\infty)$, the zeta function $\zeta^-(X(M),z)$ is obtained.\hfill$\Box$

\bigskip

The zero points of $\zeta^+(X(M),z)$ and $\zeta^-(X(M),z)$ are $1/(2M+1)$ and $1/(M+2)$, respectively. Therefore $h^+(X(M))$ and $h^-(X(M))$ also follow these values immediately. The number of periodic points is

\begin{pr}
\begin{eqnarray*}
\lefteqn{p_n(X(M))}\\
&&=\left\{ 
\begin{array}[c]{lll}
(2M+1)^n+(M+2)^n-\sum_{i=0}^{n/2}{n \choose i}2^i(M^i+M^{n-i})+{n \choose \frac{n}2}(2M)^\frac{n}2&\mbox{if}\ n\ \mbox{is\ even}&\\
&&\\
(2M+1)^n+(M+2)^n-2\sum_{i=0}^{(n-1)/2}{n \choose i}2^i(M^i+M^{n-i})&\mbox{if}\ n\ \mbox{is\ odd}.&
\end{array}\right.
\end{eqnarray*}
\end{pr}

\noindent{\bf Proof}\quad For $n\geq1$,
$p_n(X(M))=p_n(E_+^\infty)+p_n(E_-^\infty)+p_n(R_1^\infty)+p_n(L_1^\infty)-p_n(E^\infty).$\\
Since $p_n(X_{G'_1})=(2M+1)^n$ and $p_n(X_{G''_1})=(M+2)^n$, $n\geq1$ from (7) and (8), 
$$p_n(X(M))=(2M+1)^n+(M+2)^n-\{p_n(\widetilde E_-^\infty)+p_n(\widetilde L_1^\infty)\}-\{p_n(\widetilde E_+^\infty)+p_n(R_1^\infty)\}+p_n(E^\infty).$$
Furthermore, 
$$p_n(E^\infty)=
\left\{ 
\begin{array}[c]{ll}
{n \choose \frac{n}2}(2M)^\frac{n}2&\quad\mbox{if}\ n\ \mbox{is\ even},\\
0&\quad\mbox{if}\ n\ \mbox{is\ odd}.
\end{array}\right.
$$
\begin{eqnarray*}
p_n(\widetilde E_-^\infty)+p_n(\widetilde L_1^\infty)
&=&\sharp\left\{x\in P_n(X(M)):\begin{array}{l}
\mbox{a multiplier of}\ x\ \mbox{is neutral or}\\
\mbox{a positive power of}\ \lambda\end{array}\right\}\\
&=&\textstyle{\sum_{i=0}^{[n/2]}{n \choose i}(2M)^i}
\end{eqnarray*}
and
\begin{eqnarray*}
\lefteqn{p_n(\widetilde E_+^\infty)+p_n(\widetilde R_1^\infty)}\\
&&=\sharp\left\{x\in P_n(X(M)):\begin{array}{l}
\mbox{a multiplier of}\ x\ \mbox{is neutral or}\\
\mbox{a positive power of}\ \rho_1,\rho_2,\cdots,\rho_M\end{array}\right\}\\
&&=\textstyle{\sum_{i=0}^{[n/2]}{n \choose i}2^iM^{n-i}}.
\end{eqnarray*}
(See the proof of Theorem 2.5 in [HI],[I].) Thus the proof is complete. \hfill$\Box$



\end{document}